\documentclass[12pt,leqno]{amsart}

\usepackage{amscd}
\usepackage{amsmath}
\usepackage{amsfonts}
\usepackage{amssymb}
\usepackage{url}

\title{Noether's problem for $p$-groups with an abelian subgroup of index $p$}
\author{Ivo M. Michailov}
\address{Faculty of Mathematics and Informatics, Constantin Preslavski University, Universitetska str. 115, 9700 Shumen, Bulgaria}
\email{ivo\_michailov@yahoo.com}
\date{\today}
\thanks{This work is partially supported by a project No RD-08-241/12.03.2013 of Shumen University}
\begin{document}
\baselineskip 20pt
\begin{abstract}
Let $K$ be a field and $G$ be a finite group. Let $G$ act on the
rational function field $K(x(g):g\in G)$ by $K$-automorphisms
defined by $g\cdot x(h)=x(gh)$ for any $g,h\in G$. Denote by $K(G)$
the fixed field $K(x(g):g\in G)^G$. Noether's problem then asks
whether $K(G)$ is rational over $K$. Let $p$ be an odd prime and let
$G$ be a $p$-group of exponent $p^e$. Assume also that {\rm (i)}
char $K = p>0$, or {\rm (ii)} char $K \ne p$ and $K$ contains a
primitive $p^e$-th root of unity. In this paper we prove that $K(G)$
is rational over $K$ for the following two types of groups: {\rm
(1)} $G$ is a finite $p$-group with an abelian normal subgroup $H$
of index $p$, such that $H$ is a direct product of normal subgroups
of $G$ of the type $C_{p^b}\times (C_p)^c$ for some $b,c:1\leq
b,0\leq c$; {\rm (2)} $G$ is any group of order $p^5$ from the
isoclinic families with numbers $1,2,3,4,8$ and $9$.
\end{abstract}

\maketitle
\newcommand{\Gal}{{\rm Gal}}
\newcommand{\Ker}{{\rm Ker}}
\newcommand{\GL}{{\rm GL}}
\newcommand{\Br}{{\rm Br}}
\newcommand{\lcm}{{\rm lcm}}
\newcommand{\ord}{{\rm ord}}
\numberwithin{equation}{section}
\renewcommand{\thefootnote}{\fnsymbol{footnote}}

%--------------------------------------Section 1-------------------------------------------------
\section{Introduction}
\label{1}

Let $K$ be a field and $G$ be a finite group. Let $G$ act on the
rational function field $K(x(g):g\in G)$ by $K$-automorphisms
defined by $g\cdot x(h)=x(gh)$ for any $g,h\in G$. Denote by $K(G)$
the fixed field $K(x(g):g\in G)^G$. {\it Noether's problem} then
asks whether $K(G)$ is rational ($=$ purely transcendental) over
$K$. It is related to the inverse Galois problem, to the existence
of generic $G$-Galois extensions over $k$, and to the existence of
versal $G$-torsors over $k$-rational field extensions (see
\cite{Sw,Sa1} and \cite[33.1, p.86]{GMS}).

The following well-known theorem gives a positive answer to the
Noether's problem for abelian groups over a field $K$ which contains
enough roots of unity.

\newtheorem{t1.1}{Theorem}[section]
\begin{t1.1}\label{t1.1}
{\rm (Fischer }\cite[Theorem 6.1]{Sw}{\rm )} Let $G$ be a finite
abelian group of exponent $e$. Assume that {\rm (i)} either char $K
= 0$ or char $K > 0$ with char $K\nmid e$, and {\rm (ii)} $K$
contains a primitive $e$-th root of unity. Then $K(G)$ is rational
over $K$.
\end{t1.1}

Swan's paper \cite{Sw} also gives a survey of many results related
to the Noether's problem for abelian groups. In the same time, just
a handful of results about Noether's problem are obtained when the
groups are non-abelian.

We are going to list several results obtained recently by Kang:

\newtheorem{t1.2}[t1.1]{Theorem}
\begin{t1.2}\label{t1.2}
{\rm (}\cite[Theorem 1.5]{Ka1}{\rm )} Let $G$ be a metacyclic
$p$-group with exponent $p^e$, and let $K$ be any field such that
{\rm (i)} char $K = p$, or {\rm (ii)} char $K \ne p$ and $K$
contains a primitive $p^e$-th root of unity. Then $K(G)$ is rational
over $K$.
\end{t1.2}

\newtheorem{t1.21}[t1.1]{Theorem}
\begin{t1.21}\label{t1.21}
{\rm (}\cite[Theorem 1.8]{Ka2}{\rm )} Let $n\geq 3$ and let $G$ be a
non-abelian $p$-group of order $p^n$ such that $G$ contains a cyclic
subgroup of index $p^2$. Assume that $K$ is any field satisfying
that either {\rm (i)} char $K = p>0$, or {\rm (ii)} char $K \ne p$
and $K$ contains a primitive $p^{n-2}$-th root of unity. Then $K(G)$
is rational over $K$.
\end{t1.21}

\newtheorem{t1.22}[t1.1]{Theorem}
\begin{t1.22}\label{t1.22}
{\rm (}\cite[Cor. 3.2]{Ka3}{\rm )} Let $K$ be a field and $G$ be a
finite group. Assume that {\rm (i)} $G$ contains an abelian normal
subgroup $H$ so that $G/H$ is cyclic of order $n$, {\rm (ii)}
$\mathbb Z[\zeta_n]$ is a unique factorization domain, and {\rm
(iii)} $\zeta_{e'}\in K$ where $e'={\rm lcm}\{\ord(\tau),\exp(H)\}$
and $\tau$ is some element of $G$ whose image generates $G/H$. If
$G\rightarrow \GL(V)$ is any finite-dimensional linear
representation of $G$ over $K$, then $K(V)^G$ is rational over $K$.
\end{t1.22}

The reader is referred to \cite{CK,HuK} for other previous results
of Noether's problem for $p$-groups. It is still an open problem
whether Theorem \ref{t1.22} could be extended for other similar
types of meta-abelian groups with a cyclic quotient. Notice the
condition that $\mathbb Z[\zeta_n]$ is a unique factorization domain
is satisfied only for $45$ integers $n$, listed in \cite[Theorem
1.5]{Ka3} (the proof is given by Masley and Montgomery \cite{MM}).

The purpose of this paper is to extend the above results for a
certain class of $p$-groups with an abelian subgroup of index $p$.
However, we should not ''over-generalize'' the above Theorems,
because Saltman and Bogomolov prove the following results.

\newtheorem{t1.3}[t1.1]{Theorem}
\begin{t1.3}\label{t1.3}
{\rm (Saltman }\cite{Sa2}{\rm )} For any prime number $p$ and for
any infinite field $K$ with char $K \ne p$ (in particular, $K$ may
be an algebraically closed field), there is a meta-abelian $p$-group
$G$ of order $p^9$ such that $K(G)$ is not rational over $K$.
\end{t1.3}

\newtheorem{t1.31}[t1.1]{Theorem}
\begin{t1.31}\label{t1.31}
{\rm (Bogomolov }\cite{Bo}{\rm )} Let $p$ be any prime number, $k$
be any algebraically closed field with char $k\ne p$. There is a
group $G$ with order $p^6$ such that $K(G)$ is not rational over
$K$.
\end{t1.31}

Moreover, Hoshi, Kang and Kunyavskii \cite[Theorem 1.12]{HKK} proved
recently that if $G$ is a group of order $p^5$ which belongs to the
isoclinism family $\Phi_{10}$, then $K(G)$ is not rational over $K$.

Let $G$ be a group of order $p^n$ for $n\geq 2$ with an abelian
subgroup $H$ of order $p^{n-1}$. Bender \cite{Be2} determined some
properties of these groups.

We introduce some notations now. The cyclic group of order $n$ we
denote by $C_n$. The subgroups $G_{(0)}=G$ and
$G_{(i)}=[G,G_{(i-1)}]$ for $i\geq 1$ are called the lower central
series of $G$. Denote $H^p=\{h^p:h\in H\}$. Assume that $H$ is
decomposed as a product of abelian groups in the following way:
$H\simeq (C_p)^k\times C_{p^{i_1}}\times
C_{p^{i_2}}\times\cdots\times C_{p^{i_t}}$ for $1<i_1\leq
i_2\leq\cdots\leq i_t$, $k\geq 1$ and $k+i_1+i_2+\cdots +i_t=n-1$.
Let $\alpha_j$ the generator of the factor $C_{p^{i_j}}$ (for $1\leq
j\leq t$). Choose arbitrary $\alpha\in G$ such that $\alpha\notin
H$. Define the groups $H_j=\langle\alpha^{-x}\alpha_j\alpha^x:0\leq
x\leq p-1\rangle$ for $1\leq j\leq t$. In the following Proposition
we find a necessary and sufficient condition for the decomposition
of $H$ as a direct product of normal subgroups of $G$ of the type
$C_{p^b}\times (C_p)^c$.

\newtheorem{prop}[t1.1]{Proposition}
\begin{prop}\label{prop}
Let $p$ be an odd prime, let $n\geq 2$ and let $G$ be a group of
order $p^n$ with an abelian subgroup $H$ of order $p^{n-1}$. Assume
that $H$ has at least one cyclic factor of order $p$.  Then $H$ is a
direct product of normal subgroups of $G$ of the type $C_{p^b}\times
(C_p)^c$ for some $b,c:1\leq b,0\leq c$, if and only if the
following two conditions are satisfied:
\begin{enumerate}
    \item  The $p$-th lower central subgroup
$G_{(p)}$ is trivial;
    \item For any $j:1\leq j\leq t$ we have $H_j\cap H^p=H_j^p$.
\end{enumerate}
\end{prop}

In the following main result we investigate Noether's problem for
these $p$-groups. The key idea to prove our result is to find a
faithful $G$-subspace $W$ of the regular representation space
$\bigoplus_{g\in G} K\cdot x(g)$ and to show that $W^G$ is rational
over $K$. The subspace $W$ is obtained as an induced representation
from $H$. We apply a method from \cite{Ka1} to find this
representation and then by various linearization methods we reduce
the problem to verifying some combinatorial identities.

\newtheorem{t1.4}[t1.1]{Theorem}
\begin{t1.4}\label{t1.4}
Let $p$ be an odd prime. Let $G$ be a group of order $p^n$ for
$n\geq 2$ with an abelian normal subgroup $H$ of order $p^{n-1}$,
and let $G$ be of exponent $p^e$. Assume that $H=H_1\times
H_2\times\cdots\times H_s$ for some $s\geq 1$ where $H_j\simeq
C_{p^{i_j}}\times (C_p)^{k_j}$ and $H_j$ is normal in $G$ for $1\leq
j\leq s, 0\leq k_j,1\leq i_1\leq i_2\leq\cdots\leq i_s$. Assume also
that {\rm (i)} char $K = p>0$, or {\rm (ii)} char $K \ne p$ and $K$
contains a primitive $p^e$-th root of unity. Then $K(G)$ is rational
over $K$.
\end{t1.4}

The rationality problem for all $p$-groups of order $\leq p^4$ was
solved by Chu and Kang.

\newtheorem{t1.6}[t1.1]{Theorem}
\begin{t1.6}\label{t1.6}
{\rm (}\cite[Theorem 1.6]{CK}{\rm )}  Let $G$ be a $p$-group of
order $\leq p^4$, and let $G$ be of exponent $p^e$. Assume that {\rm
(i)} char $K = p>0$, or {\rm (ii)} char $K \ne p$ and $K$ contains a
primitive $p^e$-th root of unity. Then $K(G)$ is rational over $K$.
\end{t1.6}

Bender classified in \cite{Be1} the groups of order $p^5$ which
contain an abelian subgroup of order $p^4$. Their number is $p+46$
(if $p-1\ne 3k$) and $p+48$ (if $p-1=3k$). By studying the
classification of all groups of order $p^5$ made by James in
\cite{Ja}, we see that the non-abelian groups with an abelian
subgroup of order $p^4$ and that are not direct products of smaller
groups belong to the isoclinic families with numbers $2,3,4$ and
$9$. In our second main result we give a positive answer for all
groups from these families. As we show in our proof, most of these
groups satisfy the conditions of Theorem \ref{t1.4}. For the
remaining groups we show in a unified way that the answer to
Noether's problem is indeed affirmative.

\newtheorem{t1.5}[t1.1]{Theorem}
\begin{t1.5}\label{t1.5}
Let $p$ be an odd prime and let $G$ be any group of order $p^5$ from
the isoclinic families (described by James \cite{Ja}) with numbers
$1,2,3,4,8$ and $9$. Assume that $G$ is of exponent $p^e$. Assume
also that {\rm (i)} char $K = p>0$, or {\rm (ii)} char $K \ne p$ and
$K$ contains a primitive $p^e$-th root of unity. Then $K(G)$ is
rational over $K$.
\end{t1.5}

Notice that Theorem \ref{t1.5} is a generalization of \cite[Theorem
4.1]{HKK}, where is shown that $K(G)$ is \emph{retract} rational
over $K$ for any group $G$ that belongs to the isoclinism family
$\Phi_i$ where $1\leq i\leq 4$ or $8\leq i\leq 9$.

We organize this paper as follows. In Section \ref{2} we recall some
preliminaries which will be used in the proof of Theorem \ref{t1.4}.
We prove Proposition \ref{prop} in Section \ref{3}. The proofs of
Theorems \ref{t1.4} and \ref{t1.5} are given in Sections \ref{4} and
\ref{5}, respectively.

%--------------------------------------Section 2-------------------------------------------------
\section{Generalities}
\label{2}

We list several results which will be used in the sequel.

\newtheorem{t2.1}{Theorem}[section]
\begin{t2.1}\label{t2.1}
{\rm (}\cite[Theorem 1]{HK}{\rm )} Let $G$ be a finite group acting
on $L(x_1,\dots,x_m)$, the rational function field of $m$ variables
over a field $L$ such that
\begin{description}
    \item [(i)] for any $\sigma\in G, \sigma(L)\subset L;$
    \item [(ii)] the restriction of the action of $G$ to $L$ is
    faithful;
    \item [(iii)] for any $\sigma\in G$,
    \begin{equation*}
\begin{pmatrix}
\sigma(x_1)\\
\vdots\\
\sigma(x_m)\\
\end{pmatrix}
=A(\sigma)\begin{pmatrix}
x_1\\
\vdots\\
x_m\\
\end{pmatrix}
+B(\sigma)
\end{equation*}
where $A(\sigma)\in\GL_m(L)$ and $B(\sigma)$ is $m\times 1$ matrix
over $L$. Then there exist $z_1,\dots,z_m\in L(x_1,\dots,x_m)$ so
that $L(x_1,\dots,x_m)^G=L^G(z_1,\dots,z_m)$ and $\sigma(z_i)=z_i$
for any $\sigma\in G$, any $1\leq i\leq m$.
\end{description}
\end{t2.1}

\newtheorem{t2.2}[t2.1]{Theorem}
\begin{t2.2}\label{t2.2}
{\rm (}\cite[Theorem 3.1]{AHK}{\rm )} Let $G$ be a finite group
acting on $L(x)$, the rational function field of one variable over a
field $L$. Assume that, for any $\sigma\in G,\sigma(L)\subset L$ and
$\sigma(x)=a_\sigma x+b_\sigma$ for any $a_\sigma,b_\sigma\in L$
with $a_\sigma\ne 0$. Then $L(x)^G=L^G(z)$ for some $z\in L[x]$.
\end{t2.2}

\newtheorem{t2.3}[t2.1]{Theorem}
\begin{t2.3}\label{t2.3}
{\rm (Kuniyoshi }\cite[Theorem 1.7]{CK}{\rm )} If $char K=p>0$ and
$G$ is a finite $p$-group, then $K(G)$ is rational over $K$.
\end{t2.3}

Finally, we give a Lemma, which can be extracted from some proofs in
\cite{Ka2,HuK}.

\newtheorem{l2.7}[t2.1]{Lemma}
\begin{l2.7}\label{l2.7}
Let $\langle\tau\rangle$ be a cyclic group of order $n>1$, acting on
$L(v_1,\dots,v_{n-1})$, the rational function field of $n-1$
variables over a field $L$ such that
\begin{eqnarray*}
\tau&:&v_1\mapsto v_2\mapsto\cdots\mapsto v_{n-1}\mapsto (v_1\cdots
v_{n-1})^{-1}\mapsto v_1.
\end{eqnarray*}
If $L$ contains a primitive $n$th root of unity $\xi$, then
$L(v_1,\dots,v_{n-1})=L(s_1,\dots,s_{n-1})$ where $\tau:s_i\mapsto
\xi^is_i$ for $1\leq i\leq n-1$.
\end{l2.7}
\begin{proof}
Define $w_0=1+v_1+v_1v_2+\cdots+v_1v_2\cdots
v_{n-1},w_1=(1/w_0)-1/n,w_{i+1}=(v_1v_2\cdots v_i/w_0)-1/n$ for
$1\leq i\leq n-1$. Thus $L(v_1,\dots,v_{n-1})=L(w_1,\dots,w_n)$ with
$w_1+w_2+\cdots+w_n=0$ and
\begin{eqnarray*}
\tau&:&w_1\mapsto w_2\mapsto\cdots\mapsto w_{n-1}\mapsto w_n\mapsto
w_1.
\end{eqnarray*}
Define $s_i=\sum_{1\leq j\leq n}\xi^{-ij}w_j$ for $1\leq i\leq n-1$.
Then $L(w_1,\dots,w_n)=L(s_1,\dots,s_{n-1})$ and $\tau:s_i\mapsto
\xi^is_i$ for $1\leq i\leq n-1$.
\end{proof}

\newtheorem{t2.4}[t2.1]{Theorem}
\begin{t2.4}\label{t2.4}
{\rm (}\cite[Theorem 1.3]{KP}{\rm )} Let $K$ be any field, and let
$H$ and $G$ be finite groups. If $K(H)$ is rational (resp. stably
rational, retract rational) over $K$, so is $K(H\times G)$ over
$K(G)$. In particular, if both $K(H)$ and $K(G)$ are rational (resp.
stably rational, retract rational) over $K$, so is $K(H\times G)$
over K.
\end{t2.4}

%--------------------------------------Section 3-------------------------------------------------
\section{Proof of Proposition \ref{prop}}
\label{3}

\emph{I. 'If' part.} Assume that the conditions (1) and (2) from the
statement are satisfied. Put $\beta_1=\alpha_1$. Since
$G_{(p)}=\{1\}$, there exist $\beta_1,\dots,\beta_k\in H$ for some
$k:1\leq k\leq p$ such that $[\beta_j,\alpha]=\beta_{j+1}$, where
$1\leq j\leq k-1$ and $\beta_k\ne 1$ is central. If
$\langle\beta_1\rangle$ is a normal subgroup of $G$, then for
$\mathcal H=(C_p)^kH_2\cdots H_t$ we have
$H\cong\langle\beta_1\rangle\times\mathcal H$. In this way, we see
that without a loss of generality, we may assume that
$\langle\beta_1\rangle$ is not normal in $G$, and in particular
$\beta_1$ is not central (i.e., $k\geq 2$).

We are going to show now that the order of $\beta_2$ is not greater
than $p$.

From $[\beta_j,\alpha]=\beta_{j+1}$ it follows the well known
formula
\begin{equation}\label{e3.1}
\alpha^{-p}\beta_1\alpha^p=\beta_1\beta_2^{\binom{p}{1}}\beta_3^{\binom{p}{2}}\cdots
\beta_p^{\binom{p}{p-1}}\beta_{p+1},
\end{equation}
where we put $\beta_{k+1}=\cdots=\beta_{p+1}=1$. Since $\alpha^p$ is
in $H$, we obtain the formula
$$\beta_2^{\binom{p}{1}}\beta_3^{\binom{p}{2}}\cdots
\beta_k^{\binom{p}{k-1}}=1.$$ Hence $(\beta_2\cdot\prod_{j\ne
2}\beta_j^{a_j})^p=1$ for some integers $a_j$. This identity clearly
is impossible if the order of $\beta_2$ is greater than $p$.

It is not hard to conclude now that $H_1\simeq C_{p^{i_1}}\times
(C_p)^{k_1}$ for some $k_1\geq 1$. Note that the elements of $H_1$
are not $p$-th powers of the elements from $H_2\cdots H_t$, since
$H_1\cap H^p=H_1^p$. Furthermore, we can adjust the generators of
$H_2,\dots,H_t$ so that $H_1\cap(H_2\cdots H_t)=\{1\}$. For example,
if we assume that $[\alpha,\alpha_1]=[\alpha,\alpha_2]$, we can
define $\alpha_2'=\alpha_2\alpha_1^{-1}$ and get
$[\alpha,\alpha_2']=1$. Define $\mathcal
H_2=\langle\alpha^{-x}\alpha_2'\alpha^x\rangle$. Clearly,
$H=(C_p)^kH_1\mathcal H_2\cdots H_t$ and $H_1\cap \mathcal
H_2=\{1\}$. With similar changes of the generators we can treat the
more general case $[\alpha^x,\alpha_1]=[\alpha^y,\alpha_2]^z$.
Proceeding by induction we will obtain a decomposition
$H=(C_p)^k\mathcal H_1\mathcal H_2\cdots \mathcal H_t$ such that
$\mathcal H_j\cap(\mathcal H_{j+1}\cdots\mathcal H_t)=\{1\}$ for any
$j$. Therefore $H=N_1\times\cdots \times N_r\times\mathcal
H_1\times\cdots\times\mathcal H_t$ where $N_1,\dots,N_r$ are normal
groups of the type $(C_p)^a$.

\emph{II. 'Only if' part.} Assume that $H=N_1\times\cdots \times
N_r\times\mathcal H_1\times\cdots\times\mathcal H_t$ where
$N_1,\dots,N_r$ are normal groups of the type $(C_p)^a$ and
$\mathcal H_1,\dots,\mathcal H_t$ are normal groups of the type
$C_{p^b}\times (C_p)^c$.

Suppose that $G_{(p)}\ne\{1\}$. Then we can assume that there exist
$\beta_1,\dots,\beta_{p+1}\in\mathcal H_1$ such that
$[\beta_j,\alpha]=\beta_{j+1}$, where $1\leq j\leq p$ and
$\beta_{p+1}\ne 1$. (We again assume that $\beta_1$ is the generator
of the factor of the type $C_{p^b}$.) From the identity \eqref{e3.1}
it follows that
$$\beta_2^{\binom{p}{1}}\beta_3^{\binom{p}{2}}\cdots
\beta_p^{\binom{p}{p-1}}\beta_{p+1}=1.$$ Hence $\beta_2$ will have
an order bigger than $p$, which is a contradiction.  Therefore,
$G_{(p)}=\{1\}$.

Now, suppose that $H_j\cap H^p\ne H_j^p$. Let us consider first the
particular case $H_j\leq \mathcal H_j$. Then there exists
$x\in\mathbb Z$ such that $[\alpha_j,\alpha^x]\notin
\langle\alpha_j^p\rangle$ and $[\alpha_j,\alpha^x]\in H_j\cap
H^p=\mathcal H_j^p$. Hence $[\alpha_j,\alpha^x]$ is of order not
less than $p^2$, so $\mathcal H_j$ can not be of the type
$C_{p^b}\times (C_p)^c$, which is a contradiction.

Finally, consider the general case when $\alpha_j$ is decomposed in
an unique way as a product of elements from some normal factors
$N_i$ and $\mathcal H_k$ ($1\leq i\leq r,1\leq k\leq t$). Then there
exists $x\in\mathbb Z$ such that $[\alpha_j,\alpha^x]=\beta^p$ for
some $\beta\in H$, where $\beta\notin H_j^p$. From the uniqueness of
the decompositions of $[\alpha_j,\alpha^x]$ and $\beta$ as products
of elements from the normal factors, it follows that we can reduce
this general case to the particular case we just considered. We are
done.

%--------------------------------------Section 4-------------------------------------------------
\section{Proof of Theorem \ref{t1.4}}
\label{4}

If char $K=p>0$, we can apply Theorem \ref{t2.3}. Therefore, we will
assume that char $K\ne p$.

Recall that $H=H_1\times H_2\times\cdots\times H_s$, where
$H_j\simeq C_{p^{i_j}}\times (C_p)^{k_j}$. Denote by $\beta_1$ the
generator of the direct factor $C_{p^{i_1}}$ and put $k=k_1$. Then
there exist $\beta_2,\dots,\beta_k\in H_1$  such that
$[\beta_j,\alpha]=\beta_{j+1}$, where $1\leq j\leq k-1$ and
$\beta_k\ne 1$ is central.

We divide the proof into several steps. We are going now to find a
faithful representation of $G$.

\emph{Step 1.} Let $V$ be a $K$-vector space whose dual space $V^*$
is defined as $V^*=\bigoplus_{g\in G}K\cdot x(g)$ where $G$ acts on
$V^*$ by $h\cdot x(g)=x(hg)$ for any $h,g\in G$. Thus
$K(V)^G=K(x(g):g\in G)^G=K(G)$.

Define $X_1,X_2,\dots,X_k\in V^*$ by
\begin{equation*}
X_j=\sum_{\ell_1,\dots,\ell_k}x\left(\prod_{m\ne
j}\beta_m^{\ell_m}\right),
\end{equation*}
for $1\leq j\leq k$. Note that $\beta_j\cdot X_i=X_i$ for $j\ne i$.
Let $\zeta_{p^{i_1}}\in K$ be a primitive $p^{i_1}$-th root of unity
and let $\zeta$ be a primitive $p$-th root of unity. Define
$Y_1,Y_2,\dots,Y_k\in V^*$ by
\begin{equation*}
Y_1=\sum_{r=0}^{p^{i_1}-1}\zeta_{p^{i_1}}^{-r}\beta_1^r\cdot X_1,~
Y_j=\sum_{r=0}^{p-1}\zeta^{-r}\beta_j^r\cdot X_j
\end{equation*}
for $2\leq j\leq k$. \footnote{If $\beta_1$ is central, i.e., $k=1$
we simply do not have $\beta_j$ and $Y_j$ for $j\geq 2$. The proof
remains valid, however.}

It follows that {\allowdisplaybreaks\begin{align*} \beta_1\ :\
&Y_1\mapsto\zeta_{p^{i_1}} Y_1,~ Y_i\mapsto Y_i,\ \text{for}\ i\ne
1,\\
 \beta_j\ :\
&Y_j\mapsto\zeta Y_j,~ Y_i\mapsto Y_i,\ \text{for}\ i\ne j\
\text{and}\ 2\leq j\leq k.
\end{align*}}
Thus $V_1=\bigoplus_{1\leq j\leq k}K\cdot Y_j$ is a representation
space of the subgroup $H_1$. In the same way we can construct a
representation space $V_j$ of the group $H_j$ for any $j:2\leq j\leq
s$. Therefore, $\bigoplus_{1\leq j\leq s}V_j$ is a representation
space of the subgroup $H$.

Define $x_{ji}=\alpha^i\cdot Y_j$ for $1\leq j\leq k,0\leq i\leq
p-1$. Recall that $[\beta_j,\alpha]=\beta_{j-1}$. Hence
$$\alpha^{-i}\beta_j\alpha^i=\beta_j\beta_{j+1}^{\binom{i}{1}}\beta_{j+2}^{\binom{i}{2}}\cdots
\beta_k^{\binom{i}{k-j}}.$$

It follows that {\allowdisplaybreaks\begin{align*}\beta_1\ :\
&x_{1i}\mapsto\zeta_{p^{i_1}} x_{1i},~ x_{ji}\mapsto
\zeta^{\binom{i}{j-1}} x_{ji},\ \text{for}\ 2\leq j\leq k\
\text{and}\ 0\leq i\leq p-1,\\
\beta_j\ :\ &x_{\ell i}\mapsto x_{\ell i},~ x_{mi}\mapsto
\zeta^{\binom{i}{m-j}} x_{mi},\ \text{for}\ 1\leq \ell\leq j-1,j\leq
m\leq k\ \text{and}\ 0\leq i\leq p-1,\\
\alpha\ :\ &x_{j0}\mapsto x_{j1}\mapsto\cdots\mapsto x_{jp-1}\mapsto
\zeta_{p^{a_j}}^{b_j}x_{j0},\ \text{for}\ 1\leq j\leq k,
\end{align*}}
where $a_j,b_j$ are some integers such that $0\leq b_j< p^{a_j}\leq
p^{i_1}$.

Clearly, $W_1=\bigoplus_{j,i}K\cdot x_{ij}\subset V^*$ is the
induced $G$-subspace obtained from $V_1$. In the same way we can
construct the induced subspaces $W_j$ obtained from $V_j$. We find
that $W=\bigoplus_{1\leq j\leq s}W_j$ is a faithful $G$-subspace of
$V^*$. Thus, by Theorem \ref{t2.1} it suffices to show that $W^G$ is
rational over $K$.

Next, we will consider the action of $G$ on $W_1$.

\emph{Step 2.} For $1\leq j\leq k$ and for $1\leq i\leq p-1$ define
$y_{ji}=x_{ji}/x_{ji-1}$. Thus $W_1=K(x_{j0},y_{ji}:1\leq j\leq
k,1\leq i\leq p-1)$ and for every $g\in G$
\begin{equation*}
g\cdot x_{j0}\in K(y_{ji}:1\leq j\leq k,1\leq i\leq p-1)\cdot
x_{j0},\ \text{for}\ 1\leq j\leq k
\end{equation*}
while the subfield $K(y_{ji}:1\leq j\leq k,1\leq i\leq p-1)$ is
invariant by the action of $G$, i.e.,
{\allowdisplaybreaks\begin{align*}\beta_1\ :\ &y_{1i}\mapsto
y_{1i},~ y_{ji}\mapsto \zeta^{\binom{i-1}{j-2}} y_{ji},\ \text{for}\
2\leq j\leq k\
\text{and}\ 1\leq i\leq p-1,\\
\beta_j\ :\ &y_{\ell i}\mapsto y_{\ell i},~ y_{mi}\mapsto
\zeta^{\binom{i-1}{m-j-1}} y_{mi},\ \text{for}\ 1\leq \ell\leq
j,j+1\leq
m\leq k\ \text{and}\ 1\leq i\leq p-1,\\
\alpha\ :\ &y_{j1}\mapsto y_{j2}\mapsto\cdots\mapsto y_{jp-1}\mapsto
\zeta_{p^{a_j}}^{b_j}(y_{j1}\cdots y_{jp-1})^{-1},\ \text{for}\
1\leq j\leq k.
\end{align*}}
From Theorem \ref{t2.2} it follows that if $K(y_{ji}:1\leq j\leq
k,1\leq i\leq p-1)^{G}$ is rational over $K$, so is
$K(x_{j0},y_{ji}:1\leq j\leq k,1\leq i\leq p-1)^{G}$ over $K$.

Since $K$ contains a primitive $p^e$-th root of unity $\zeta_{p^e}$
where $p^e$ is the exponent of $G$, $K$ contains as well a primitive
$p^{{a_j}+1}$-th root of unity, and we may replace the variables
$y_{ji}$ by $y_{ji}/\zeta_{p^{{a_j}+1}}^{b_j}$ so that we obtain a
more convenient action of $\alpha$ without changing the actions of
$\beta_j$'s. Namely we may assume that
\begin{align}\label{e4.1}
\alpha\ :\ &y_{j1}\mapsto y_{j2}\mapsto\cdots\mapsto y_{jp-1}\mapsto
(y_{j1}y_{j2}\dots y_{jp-1})^{-1}\ \text{for}\ 1\leq j\leq k.
\end{align}

Define $u_{k1}=y_{k1}^p,u_{ki}=y_{ki}/y_{ki-1}$ for $2\leq i\leq
p-1$. Then $K(y_{ji},u_{ki}:1\leq j\leq k-1,1\leq i\leq
p-1)=K(y_{ji}:1\leq j\leq k,1\leq i\leq
p-1)^{\langle\beta_{k-1}\rangle}$. From Theorem \ref{t2.2} it
follows that if $K(y_{ji},u_{ki}:1\leq j\leq k-1,2\leq i\leq
p-1)^{G}$ is rational over $K$, so is $K(y_{ji},u_{ki}:1\leq j\leq
k-1,1\leq i\leq p-1)^{G}$ over $K$. We have the following actions
{\allowdisplaybreaks\begin{align*}\beta_j\ :\ &u_{ki}\mapsto
\zeta^{\binom{i-2}{k-j-2}} u_{ki},\
\text{for}\ 2\leq i\leq p-1\ \text{and}\ 1\leq j\leq k-2,\\
\alpha\ :\ &u_{k2}\mapsto u_{k3}\mapsto\cdots\mapsto u_{kp-1}\mapsto
(u_{k1}u_{k2}^{p-1}u_{k3}^{p-2}\cdots u_{kp-1}^2)^{-1}\mapsto
u_{k1}u_{k2}^{p-2}u_{k3}^{p-3}\cdots u_{kp-2}^2u_{kp-1}.
\end{align*}}
For $2\leq i\leq p-1$ define
$$v_{ki}=u_{ki}y_{k-1i}^{-1}y_{k-2i}y_{k-3i}^{-1}\cdots
y_{4i}^{(-1)^k}y_{3i}^{(-1)^{k+1}}y_{2i}^{(-1)^{k+2}},$$ and put
$v_{k1}=u_{k1}$.

With the aid of the well known property
$\binom{n}{m}-\binom{n-1}{m}=\binom{n-1}{m-1}$, it is not hard to
verify the following identity {\allowdisplaybreaks\begin{align*} \
&\binom{i-2}{k-j-2}-\binom{i-1}{k-j-2}+\binom{i-1}{k-j-3}-\binom{i-1}{k-j-4}+\cdots\\
&\cdots+(-1)^{k-j-1}\binom{i-1}{2}+(-1)^{k-j}\binom{i-2}{1}=0.
\end{align*}}
It follows that {\allowdisplaybreaks\begin{align*}\beta_j\ :\
&v_{ki}\mapsto v_{ki},\
\text{for}\ 1\leq i\leq p-1\ \text{and}\ 1\leq j\leq k-2,\\
\alpha\ :\ &v_{k2}\mapsto v_{k3}\mapsto\cdots\mapsto v_{kp-1}\mapsto
A_k\cdot(v_{k1}v_{k2}^{p-1}v_{k3}^{p-2}\cdots v_{kp-1}^2)^{-1},
\end{align*}}
where $A_k$ is some monomial in $y_{ji}$ for $2\leq j\leq k-1,1\leq
i\leq p-1$.

It is obvious that we can proceed in the same way defining elements
$v_{k-1i},\dots,v_{2i}$ such that $\beta_j$ acts trivially on all
$v_{mi}$'s and the action of $\alpha$ is given by
\begin{align}\label{e4.2}
\alpha\ :\ &v_{m1}\mapsto v_{m1}v_{m2}^p,~ v_{m2}\mapsto
v_{m3}\mapsto\cdots\mapsto v_{mp-1}\mapsto
A_m\cdot(v_{m1}v_{m2}^{p-1}v_{m3}^{p-2}\cdots v_{mp-1}^2)^{-1},
\end{align}
where $A_m$ is some monomial in $v_{2i},\dots,v_{m-1i}$ for $3\leq
m\leq k$ and $A_2=1$. In this way we obtain that
$K(y_{1i},v_{ji})=K(y_{ji})^{H_1}$.

We will ''linearize'' the action \eqref{e4.2} applying repeatedly
Kang's argument from \cite[Case 5, Step II]{Ka2}. (Note that the
linearization of $\alpha$ on $y_{1i}$'s follows from Lemma
\ref{l2.7}.)

\emph{Step 3.} We write the additive version of the multiplication
action of $\alpha$ in formula \eqref{e4.1}, i.e., consider the
$\mathbb Z[\pi]$-module $M=\bigoplus_{1\leq m\leq k}(\oplus_{1\leq
i\leq p-1}\mathbb Z\cdot v_{mi})$ corresponding to \eqref{e4.2},
where $\pi=\langle\alpha\rangle$. Denote the submodules
$M_j=\bigoplus_{1\leq m\leq j}(\oplus_{1\leq i\leq p-1}\mathbb
Z\cdot v_{mi})$ for $1\leq j\leq k$. Thus $\alpha$ has the following
additive action {\allowdisplaybreaks\begin{align*}
\alpha\ :\ &v_{j1}\mapsto v_{j1}+pv_{j2},~\\
 &v_{j2}\mapsto
v_{j3}\mapsto\cdots\mapsto v_{jp-1}\mapsto
A_j-v_{j1}-(p-1)v_{j2}-(p-2)v_{j3}-\cdots -2v_{jp-1},
\end{align*}}
where $A_j\in M_{j-1}$.

By Lemma \ref{l2.7}, $M_1$ is isomorphic to the $\mathbb
Z[\pi]$-module $N=\oplus_{1\leq i\leq p-1}\mathbb Z\cdot s_i$ where
$s_1=v_{12},s_i=\alpha^{i-1}\cdot v_{12}$ for $2\leq i\leq p-1$, and
\begin{align*}
\alpha\ :\  &s_1\mapsto s_2\mapsto\cdots\mapsto s_{p-1}\mapsto
-s_1-s_2-\cdots-s_{p-1}\mapsto s_1.
\end{align*}

Let $\Phi_p(T)\in\mathbb Z[T]$ be the $p$-th cyclotomic polynomial.
Since $\mathbb Z[\pi]\simeq\mathbb Z[T]/T^p-1$, we find that
$\mathbb Z[\pi]/\Phi_p(\alpha)\simeq \mathbb Z[T]/\Phi_p(T)\simeq
\mathbb Z[\omega]$, the ring of $p$-th cyclotomic integer. As
$\Phi_p(\alpha)\cdot x=0$ for any $x\in N$, the $\mathbb
Z[\pi]$-module $N$ can be regarded as a $\mathbb Z[\omega]$-module
through the morphism $\mathbb Z[\pi]\to\mathbb
Z[\pi]/\Phi_p(\alpha)$. When $N$ is regarded as a $\mathbb
Z[\omega]$-module, $N\simeq\mathbb Z[\omega]$ the rank-one free
$\mathbb Z[\omega]$-module.

We claim that $M$ itself can be regarded as a $\mathbb
Z[\omega]$-module, i.e., $\Phi_p(\alpha)\cdot M=0$.

Return to the multiplicative notations in Step 2. Note that all
$v_{ji}$'s are monomials in $y_{ji}$'s. The action of $\alpha$ on
$y_{ji}$ given in formula \eqref{e4.1} satisfies the relation
$\prod_{0\leq m\leq p-1}\alpha^m(y_{ji})=1$ for any $1\leq j\leq
k,1\leq i\leq p-1$. Using the additive notations, we get
$\Phi_p(\alpha)\cdot y_{ji}=0$. Hence $\Phi_p(\alpha)\cdot M=0$.

Define $M'=M/M_{k-1}$. It follows that we have a short exact
sequence of $\mathbb Z[\pi]$-modules
\begin{equation}\label{e4.3}
0\to M_{k-1}\to M\to M'\to 0.
\end{equation}
Since $M$ is a $\mathbb Z[\omega]$-module, \eqref{e4.3} is a short
exact sequence of $\mathbb Z[\omega]$-modules. Proceeding by
induction, we obtain that $M$ is a direct sum of free $\mathbb
Z[\omega]$-modules isomorphic to $N$. Therefore,
$M\simeq\oplus_{1\leq j\leq k}N_j$, where $N_j\simeq N$ is a free
$\mathbb Z[\omega]$-module, and so a $\mathbb Z[\pi]$-module also
(for $1\leq j\leq k$).

Finally, we interpret the additive version of $M\simeq\oplus_{1\leq
j\leq k}N_j\simeq N^k$ it terms of the multiplicative version as
follows: There exist $w_{ji}$ that are monomials in $v_{ji}$ for
$1\leq j\leq k,1\leq i\leq p-1$ such that $K(w_{ji})=K(v_{ji})$ and
$\alpha$ acts as
\begin{align*}
\alpha\ :\ &w_{j1}\mapsto w_{j2}\mapsto\cdots\mapsto w_{jp-1}\mapsto
(w_{j1}w_{j2}\dots w_{jp-1})^{-1}\ \text{for}\ 1\leq j\leq k.
\end{align*}
According to Lemma \ref{l2.7}, the above action can be linearized.
Since $H\simeq H_1\times\cdots \times H_s$ for some normal subgroups
$H_j$ of $G$, we obtain that $W^H$ is a $K$-free compositum of
fields having a linear action of $\alpha$. Therefore, $W^G$ is
rational over $K$. We are done.

%--------------------------------------Section 5-------------------------------------------------
\section{Proof of Theorem \ref{t1.5}}
\label{5}

All groups from family 1 are abelian, so we may apply Theorem
\ref{t1.1}.

For the groups that are direct product of smaller groups (contained
in families 2 and 3) we may apply Theorems \ref{t1.6} and
\ref{t2.4}.

The groups $\Phi_2(41),\Phi_2(32)a_1,\Phi_2(32)a_2$ and $\Phi_8(32)$
are metacyclic, so we may apply Theorem \ref{t1.2}.

The group $\Phi_2(311)b$ contains a normal subgroup
$\langle\alpha_1,\gamma\rangle\simeq C_p\times C_{p^3}$. The group
$\Phi_2(311)c$ contains a normal subgroup
$\langle\alpha_2,\alpha\rangle\simeq C_p\times C_{p^3}$. The group
$\Phi_2(221)c$ contains a normal subgroup
$\langle\gamma,\alpha_1,\alpha^p\rangle\simeq C_{p^2}\times
(C_p)^2$. The group $\Phi_2(221)d$ contains a normal subgroup
$\langle\alpha_1,\alpha_2,\alpha^p\rangle\simeq C_{p^2}\times
(C_p)^2$. For all these groups we may apply Theorem \ref{t1.4}.

The group $\Phi_3(311)a$ contains a normal subgroup
$\langle\alpha^p,\alpha_1,\alpha_2\rangle\simeq C_{p^2}\times
(C_p)^2$. The group $\Phi_3(311)b_r$ contains a normal subgroup
$\langle\alpha_1,\alpha_2\rangle\simeq C_{p^3}\times C_p$. The
groups $\Phi_3(221)a$ and $\Phi_3(221)b_r$ contain a normal subgroup
$\langle\alpha_1,\alpha_2,\alpha^p\rangle\simeq C_{p^2}\times
(C_p)^2$. The group $\Phi_3(2111)c$ contains a normal subgroup
$\langle\gamma,\alpha_1,\alpha_2\rangle\simeq C_{p^2}\times
(C_p)^2$. The group $\Phi_3(2111)d$ contains a normal subgroup
$\langle\alpha_1,\alpha_2,\alpha_3,\alpha^p\rangle\simeq (C_p)^4$.
The group $\Phi_3(2111)e$ contains a normal subgroup
$\langle\alpha_1,\alpha_2,\alpha_3\rangle\simeq C_{p^2}\times
(C_p)^2$. We may again apply Theorem \ref{t1.4}.

In the same way, we see that Theorem \ref{t1.4} may be applied for
the groups $\Phi_4(221)a$, $\Phi_4(221)b$, $\Phi_4(2111)a$,
$\Phi_4(2111)b,\Phi_4(2111)c$ and $\Phi_4(1^5)$, because each group
contains a normal subgroup $H\simeq (C_p)^4$ or $C_{p^2}\times
(C_p)^2$. Similarly, for $p\geq 5$ each group from family $9$
contains a normal subgroup
$\langle\alpha_1,\alpha_2,\alpha_3,\alpha_4\rangle\simeq (C_p)^4$.
For $p=3$ each group from family $9$ contains a normal subgroup
$\langle\alpha_1,\alpha_2,\alpha_3,\alpha_4\rangle\simeq C_9\times
C_9$.

Thus it remains to consider the groups $\Phi_4(221)c$,
$\Phi_4(221)d_r$, $\Phi_4(221)e$, $\Phi_4(221)f_0$, $\Phi_4(221)f_r$
(for any odd $p$) and the groups from family $9$ for $p=3$.

Let $G$ be any of these groups and let
$H=\langle\alpha_1,\alpha_2\rangle\simeq C_{p^2}\times C_{p^2}$. Let
$V$ be a $K$-vector space whose dual space $V^*$ is defined as
$V^*=\bigoplus_{g\in G}K\cdot x(g)$ where $G$ acts on $V^*$ by
$h\cdot x(g)=x(hg)$ for any $h,g\in G$. Thus $K(V)^G=K(x(g):g\in
G)^G=K(G)$. Let $\zeta_{p^2}$ be a primitive $p^2$-th root of unity,
and let $\zeta=\zeta_{p^2}^p$. We may find $Y_1,Y_2\in V^*$ such
that {\allowdisplaybreaks\begin{align*} \alpha_1\ :\
&Y_1\mapsto\zeta_{p^2} Y_1,~ Y_2\mapsto Y_2,\\
 \alpha_2\ :\
&Y_1\mapsto Y_2,~ Y_2\mapsto\zeta_{p^2} Y_2.
\end{align*}}
Thus $V_1=\bigoplus_{1\leq j\leq 2}K\cdot Y_j$ is a representation
space of the subgroup $H$. Define $x_i=\alpha^i\cdot
Y_1,y_i=\alpha^i\cdot Y_2$ for $0\leq i\leq p-1$. Note that for the
groups from the family $\Phi_4$ we have the relations
$\alpha_j\alpha^i=\alpha^i\alpha_j\beta_j^i$ for $0\leq i\leq
p-1,1\leq j\leq 2$. From now on we are going to consider each group
individually.

\emph{Case I.} $G=\Phi_4(221)d_r$. Let $\ell$ be an integer such
that $k\ell\equiv 1\pmod{p}$. Hence $\beta_1=\alpha_1^{\ell p}$. We
have the actions {\allowdisplaybreaks\begin{align*}\alpha_1\ :\
&x_i\mapsto\zeta_{p^2}^{1+i\ell p} x_i,~ y_i\mapsto y_i,\\
\alpha_2\ :\
&x_i\mapsto x_i,~ y_i\mapsto\zeta_{p^2}^{1+ip} y_i,\\
\alpha\ :\ &x_0\mapsto x_1\mapsto\cdots\mapsto x_{p-1}\mapsto x_0,\\
&y_0\mapsto y_1\mapsto\cdots\mapsto y_{p-1}\mapsto y_0,
\end{align*}}
where $0\leq i\leq p-1$. We find that $W=(\bigoplus_i K\cdot
x_i)\bigoplus(\bigoplus_i K\cdot y_i)\subset V^*$ is a faithful
$G$-subspace of $V^*$. Thus, by Theorem \ref{t2.1} it suffices to
show that $W^G$ is rational over $K$.

For $1\leq i\leq p-1$ define $u_i=x_i/x_{i-1},v_i=y_i/y_{i-1}$ and
$\zeta=\zeta_{p^2}^p$. Thus $W=K(x_0,y_0,u_i,v_i:1\leq i\leq p-1)$
and for every $g\in G$
\begin{equation*}
g\cdot x_0\in K(u_i:1\leq i\leq p-1)\cdot x_0,\ \text{and}\ g\cdot
y_0\in K(v_i:1\leq i\leq p-1)\cdot y_0,
\end{equation*}
while the subfield $K(u_i,v_i:1\leq i\leq p-1)$ is invariant by the
action of $G$, i.e., {\allowdisplaybreaks\begin{align*}\alpha_1\ :\
&u_i\mapsto\zeta^\ell u_i,~ v_i\mapsto v_i,\\
\alpha_2\ :\
&u_i\mapsto u_i,~ v_i\mapsto\zeta v_i,\\
\alpha\ :\ &u_1\mapsto u_2\mapsto\cdots\mapsto u_{p-1}\mapsto (u_1u_2\cdots u_{p-1})^{-1},\\
&v_1\mapsto v_2\mapsto\cdots\mapsto v_{p-1}\mapsto (v_1v_2\cdots
v_{p-1})^{-1}.
\end{align*}}
From Theorem \ref{t2.2} it follows that if $K(u_i,v_i:1\leq i\leq
p-1)^{G}$ is rational over $K$, so is $K(x_0,y_0,u_i,v_i:1\leq i\leq
p-1)^{G}$ over $K$.

Define $U_1=u_1^p,U_i=u_i/u_{i-1}$ and $V_1=v_1^p,V_i=v_i/v_{i-1}$
for $2\leq i\leq p-1$. Then $K(U_i,V_i:1\leq i\leq
p-1)=K(u_i,v_i:1\leq i\leq p-1)^{\langle\alpha_1,\alpha_2\rangle}$
and the action of $\alpha$ on $U_i$ and $V_i$ is
{\allowdisplaybreaks
\begin{align*}
\alpha\ :\ &U_1\mapsto U_1U_2^p,\\
&U_2\mapsto U_3\mapsto\cdots\mapsto U_{p-1}\mapsto
(U_1U_2^{p-1}U_3^{p-2}\cdots U_{p-1}^2)^{-1}\mapsto\\
&\mapsto U_1U_2^{p-2}U_3^{p-3}\cdots U_{p-2}^2U_{p-1}\mapsto U_2,\\
&V_1\mapsto V_1V_2^p,\\
&V_2\mapsto V_3\mapsto\cdots\mapsto V_{p-1}\mapsto
(V_1V_2^{p-1}V_3^{p-2}\cdots V_{p-1}^2)^{-1}\mapsto\\
&\mapsto V_1V_2^{p-2}V_3^{p-3}\cdots V_{p-2}^2V_{p-1}\mapsto V_2.
\end{align*}}
Define $W_1=U_2,W_i=\alpha^i\cdot U_2$ and
$Z_1=V_2,Z_i=\alpha^i\cdot V_2$ for $2\leq i\leq p-1$. Now the
action of $\alpha$ is {\allowdisplaybreaks
\begin{align*}
\alpha\ :\ &W_1\mapsto W_2\mapsto\cdots\mapsto W_{p-1}\mapsto
(W_1W_2\cdots W_{p-1})^{-1},\\
&Z_1\mapsto Z_2\mapsto\cdots\mapsto Z_{p-1}\mapsto (Z_1Z_2\cdots
Z_{p-1})^{-1}.
\end{align*}}
Since $U_1=(W_{p-1}W_1^{p-1}W_2^{p-2}\cdots W_{p-2}^2)^{-1}$ and
$V_1=(Z_{p-1}Z_1^{p-1}Z_2^{p-2}\cdots Z_{p-2}^2)^{-1}$, we get that
$K(U_1,\dots,U_{p-1},V_1,\dots,V_{p-1})=K(W_1,\dots,W_{p-1},Z_1,\dots,Z_{p-1})$.
From Lemma \ref{l2.7} it follows that the action of $\alpha$ on
$K(W_1,\dots,W_{p-1},Z_1,\dots,Z_{p-1})$ can be linearized. It
remains to apply Theorem \ref{t1.1}.

\emph{Case II.} $G=\Phi_4(221)c$. The actions of $\alpha_1$ and
$\alpha_2$ on $u_i$ and $v_i$ (we keep the notations from Case I)
are {\allowdisplaybreaks\begin{align*}\alpha_1\ :\
&u_i\mapsto\zeta u_i,~ v_i\mapsto v_i,\\
\alpha_2\ :\ &u_i\mapsto u_i,~ v_i\mapsto\zeta v_i.
\end{align*}}
These actions are the same as in Case I, so we may apply the same
proof.

\emph{Case III.} $G=\Phi_4(221)f_0$. Let $\mu$ be an integer such
that $\mu\nu\equiv 1\pmod{p}$. Then $\beta_1=\alpha_2^{\mu p}$. The
actions of $\alpha_1$ and $\alpha_2$ on $u_i$ and $v_i$ are
{\allowdisplaybreaks\begin{align*}\alpha_1\ :\
&u_i\mapsto u_i,~ v_i\mapsto\zeta^\mu v_i,\\
\alpha_2\ :\ &u_i\mapsto\zeta u_i,~ v_i\mapsto v_i.
\end{align*}}
These actions are almost the same as in Case I, so we may apply the
same proof.

\emph{Case IV.} $G=\Phi_4(221)e$. Let $s$ be an integer such that
$(-1/4)s\equiv 1\pmod{p}$. Then $\beta_2=\alpha_1^{sp}$ and
$\beta_1=\alpha_2^p\alpha_1^{-sp}$. The actions of $\alpha_1$ and
$\alpha_2$ on $u_i$ and $v_i$ are
{\allowdisplaybreaks\begin{align*}\alpha_1\ :\
&u_i\mapsto\zeta^{-s} u_i,~ v_i\mapsto\zeta v_i,\\
\alpha_2\ :\ &u_i\mapsto\zeta^s u_i,~ v_i\mapsto v_i.
\end{align*}}
For $1\leq i\leq p-1$ define $w_i=u_iv_i^s$. We have now
{\allowdisplaybreaks\begin{align*}\alpha_1\ :\
&w_i\mapsto w_i,~ v_i\mapsto\zeta v_i,\\
\alpha_2\ :\ &w_i\mapsto\zeta^s w_i,~ v_i\mapsto v_i.
\end{align*}}
We may apply again the proof of Case I.

\emph{Case V.} $G=\Phi_4(221)f_r$. The proof  is the same as Case
IV.

\emph{Case VI.} $G=\Phi_9(1^5)$ for $p=3$. Calculations show that we
have the following relations:
$\alpha_1^3=\alpha_3^{-1}\alpha_4,\alpha_2^3=\alpha_4^{-1},\alpha_3^3=\alpha_4^3=1,
\alpha_1\alpha=\alpha\alpha_1\alpha_2,\alpha_1\alpha^2=\alpha^2\alpha_1^{-2}\alpha_2^{-1},
\alpha_2\alpha=\alpha\alpha_1^{-3}\alpha_2^{-2},\alpha_2\alpha^2=\alpha^2\alpha_2\alpha_1^3$.
Hence we obtain the actions
{\allowdisplaybreaks\begin{align*}\alpha_1\ :\
&x_0\mapsto\zeta_9 x_0,~ x_1\mapsto\zeta_9 x_1,~ x_2\mapsto\zeta_9^{-2} x_2,~ y_0\mapsto y_0,~ y_1\mapsto\zeta_9 y_1,~ y_2\mapsto\zeta_9^{-1} y_2,\\
\alpha_2\ :\
&x_0\mapsto x_0,~ x_1\mapsto\zeta^{-1} x_1,~ x_2\mapsto\zeta x_2,~ y_0\mapsto\zeta_9 y_0,~ y_1\mapsto\zeta_9^{-2} y_1,~ y_2\mapsto\zeta_9 y_2,\\
\alpha\ :\ &x_0\mapsto x_1\mapsto x_2\mapsto x_0,\\
&y_0\mapsto y_1\mapsto y_2\mapsto y_0,
\end{align*}}
where $\zeta_9$ is a primitive $9$-th root of unity and
$\zeta=\zeta_9^3$ is a primitive $3$-th root of unity. We find that
$W=(\bigoplus_i K\cdot x_i)\bigoplus(\bigoplus_i K\cdot y_i)\subset
V^*$ is a faithful $G$-subspace of $V^*$. Thus, by Theorem
\ref{t2.1} it suffices to show that $W^G$ is rational over $K$.

For $1\leq i\leq 2$ define $u_i=x_i/x_{i-1},v_i=y_i/y_{i-1}$. Thus
$W=K(x_0,y_0,u_i,v_i:1\leq i\leq 2)$ and for every $g\in G$
\begin{equation*}
g\cdot x_0\in K(u_i:1\leq i\leq 2)\cdot x_0,\ \text{and}\ g\cdot
y_0\in K(v_i:1\leq i\leq 2)\cdot y_0,
\end{equation*}
while the subfield $K(u_i,v_i:1\leq i\leq 2)$ is invariant by the
action of $G$, i.e., {\allowdisplaybreaks\begin{align*}\alpha_1\ :\
&u_1\mapsto u_1,~ u_2\mapsto\zeta^{-1} u_2,~ v_1\mapsto\zeta_9 v_1,~ v_2\mapsto\zeta_9^{-2} v_2,\\
\alpha_2\ :\
&u_1\mapsto\zeta^{-1} u_1,~ u_2\mapsto\zeta^{-1} u_2,~ v_1\mapsto\zeta^{-1} v_1,~ v_2\mapsto\zeta v_2,\\
\alpha\ :\ &u_1\mapsto u_2\mapsto (u_1u_2)^{-1},\\
&v_1\mapsto v_2\mapsto (v_1v_2)^{-1}.
\end{align*}}
From Theorem \ref{t2.2} it follows that if $K(u_i,v_i:1\leq i\leq
2)^{G}$ is rational over $K$, so is $K(x_0,y_0,u_i,v_i:1\leq i\leq
2)^{G}$ over $K$.

Define $w_1=v_1^3,w_2=v_2/v_1$. Then
$K(u_1,u_2,w_1,w_2)=K(u_1,u_2,v_1,v_2)^{\langle\alpha_1^3\rangle}$,
and the actions of $\alpha_1,\alpha_2$ and $\alpha$ on
$K(u_1,u_2,w_1,w_2)$ are {\allowdisplaybreaks\begin{align*}\alpha_1\
:\
&u_1\mapsto u_1,~ u_2\mapsto\zeta^{-1} u_2,~ w_1\mapsto\zeta w_1,~ w_2\mapsto\zeta^{-1} w_2,\\
\alpha_2\ :\
&u_1\mapsto\zeta^{-1} u_1,~ u_2\mapsto\zeta^{-1} u_2,~ w_1\mapsto w_1,~ w_2\mapsto\zeta^{-1} w_2,\\
\alpha\ :\ &u_1\mapsto u_2\mapsto (u_1u_2)^{-1},\\
&w_1\mapsto w_2^3w_1,~ w_2\mapsto (w_1w_2^2)^{-1}.
\end{align*}}
Define $V_1=w_2,V_2=(w_1w_2^2)^{-1}$. Then
$K(u_1,u_2,V_1,V_2)=K(u_1,u_2,w_1,w_2)$, and the actions of
$\alpha_1,\alpha_2$ and $\alpha$ on $K(u_1,u_2,V_1,V_2)$ are
{\allowdisplaybreaks\begin{align*}\alpha_1\ :\
&u_1\mapsto u_1,~ u_2\mapsto\zeta^{-1} u_2,~ V_1\mapsto\zeta^{-1} V_1,~ V_2\mapsto\zeta V_2,\\
\alpha_2\ :\
&u_1\mapsto\zeta^{-1} u_1,~ u_2\mapsto\zeta^{-1} u_2,~ V_1\mapsto\zeta^{-1} V_1,~ V_2\mapsto\zeta^{-1} V_2,\\
\alpha\ :\ &u_1\mapsto u_2\mapsto (u_1u_2)^{-1},\\
&V_1\mapsto V_2\mapsto (V_1V_2)^{-1}.
\end{align*}}
Define $U_1=u_1^3,U_2=u_2/u_1,W_1=V_1/u_1,W_2=V_2/u_2$. Then
$K(U_1,U_2,W_1,W_2)=K(u_1,u_2,V_1,V_2)^{\langle\alpha_2\rangle}$,
and the actions of $\alpha_1$ and $\alpha$ on $K(U_1,U_2,W_1,W_2)$
are {\allowdisplaybreaks\begin{align*}\alpha_1\ :\
&U_1\mapsto U_1,~ U_2\mapsto\zeta^{-1} U_2,~ W_1\mapsto\zeta^{-1} W_1,~ W_2\mapsto\zeta^{-1} W_2,\\
\alpha\ :\ &U_1\mapsto U_2^3U_1,~ U_2\mapsto (U_1U_2^2)^{-1},\\
&W_1\mapsto W_2\mapsto (W_1W_2)^{-1}.
\end{align*}}
Define $\tilde u_1=U_2,\tilde u_2=(U_1U_2^2)^{-1},\tilde
v_1=W_1/\tilde u_1,\tilde v_2=W_2/\tilde u_2$. Then $K(\tilde
u_1,\tilde u_2,\tilde v_1,\tilde v_2)=K(U_1,U_2,W_1,W_2)$, and the
actions of $\alpha_1$ and $\alpha$ on $K(\tilde u_1,\tilde
u_2,\tilde v_1,\tilde v_2)$ are
{\allowdisplaybreaks\begin{align*}\alpha_1\ :\
&\tilde u_1\mapsto\zeta^{-1} \tilde u_1,~ \tilde u_2\mapsto\zeta^{-1} \tilde u_2,~ \tilde v_1\mapsto\tilde v_1,~ \tilde v_2\mapsto\tilde v_2,\\
\alpha\ :\ &\tilde u_1\mapsto \tilde u_2\mapsto (\tilde u_1\tilde u_2)^{-1},\\
&\tilde v_1\mapsto \tilde v_2\mapsto (\tilde v_1\tilde v_2)^{-1}.
\end{align*}}
Define $\tilde U_1=\tilde u_1^3, \tilde U_2=\tilde u_2/\tilde u_1$.
Then $K(\tilde U_1,\tilde U_2,\tilde v_1,\tilde v_2)=K(\tilde
u_1,\tilde u_2,\tilde v_1,\tilde v_2)^{\langle\alpha_1\rangle}$, and
the action of $\alpha$ on $K(\tilde U_1,\tilde U_2,\tilde v_1,\tilde
v_2)$ is {\allowdisplaybreaks\begin{align*} \alpha\ :\ &\tilde
U_1\mapsto \tilde U_2^3\tilde U_1,~ \tilde U_2\mapsto (\tilde
U_1\tilde
U_2^2)^{-1},\\
&\tilde v_1\mapsto \tilde v_2\mapsto (\tilde v_1\tilde v_2)^{-1}.
\end{align*}}
Finally, define $\tilde w_1=\tilde U_2,\tilde w_2=(\tilde U_1\tilde
U_2^2)^{-1}$. Hence $\alpha(\tilde w_2)=(\tilde w_1\tilde
w_2)^{-1}$. According to Lemma \ref{l2.7} the action of $\alpha$ on
$K(\tilde w_1,\tilde w_2,\tilde v_1,\tilde v_2)$ can be linearized.

\emph{Case VII.} $G=\Phi_9(2111)a$ for $p=3$. As we have shown in
the proof of Theorem \ref{t1.4}, when $\alpha^p\in H$, we can easily
find proper substitutions so that we may assume $\alpha^p=1$. Apply
Case VI.

\emph{Case VIII.} $G=\Phi_9(2111)b_r$ for $p=3$. Calculations show
that we have the following relations:
$\alpha_1^3=\alpha_3^{-1}\alpha_4^{-1},\alpha_2^3=\alpha_4^{-1},\alpha_3^3=\alpha_4^3=1,
\alpha_1\alpha=\alpha\alpha_1\alpha_2,\alpha_1\alpha^2=\alpha^2\alpha_1^{-2}\alpha_2^{-4},
\alpha_2\alpha=\alpha\alpha_1^{-3}\alpha_2^4,\alpha_2\alpha^2=\alpha^2\alpha_2^4\alpha_1^3$.
Hence we obtain the actions
{\allowdisplaybreaks\begin{align*}\alpha_1\ :\
&x_0\mapsto\zeta_9 x_0,~ x_1\mapsto\zeta_9 x_1,~ x_2\mapsto\zeta_9^{-2} x_2,~ y_0\mapsto y_0,~ y_1\mapsto\zeta_9 y_1,~ y_2\mapsto\zeta_9^{-4} y_2,\\
\alpha_2\ :\
&x_0\mapsto x_0,~ x_1\mapsto\zeta^{-1} x_1,~ x_2\mapsto\zeta x_2,~ y_0\mapsto\zeta_9 y_0,~ y_1\mapsto\zeta_9^4 y_1,~ y_2\mapsto\zeta_9^4 y_2,\\
\alpha\ :\ &x_0\mapsto x_1\mapsto x_2\mapsto x_0,\\
&y_0\mapsto y_1\mapsto y_2\mapsto y_0.
\end{align*}}
We find that $W=(\bigoplus_i K\cdot x_i)\bigoplus(\bigoplus_i K\cdot
y_i)\subset V^*$ is a faithful $G$-subspace of $V^*$. Thus, by
Theorem \ref{t2.1} it suffices to show that $W^G$ is rational over
$K$.

For $1\leq i\leq 2$ define $u_i=x_i/x_{i-1},v_i=y_i/y_{i-1}$. Thus
$W=K(x_0,y_0,u_i,v_i:1\leq i\leq 2)$ and for every $g\in G$
\begin{equation*}
g\cdot x_0\in K(u_i:1\leq i\leq 2)\cdot x_0,\ \text{and}\ g\cdot
y_0\in K(v_i:1\leq i\leq 2)\cdot y_0,
\end{equation*}
while the subfield $K(u_i,v_i:1\leq i\leq 2)$ is invariant by the
action of $G$, i.e., {\allowdisplaybreaks\begin{align*}\alpha_1\ :\
&u_1\mapsto u_1,~ u_2\mapsto\zeta^{-1} u_2,~ v_1\mapsto\zeta_9 v_1,~ v_2\mapsto\zeta_9^4 v_2,\\
\alpha_2\ :\
&u_1\mapsto\zeta^{-1} u_1,~ u_2\mapsto\zeta^{-1} u_2,~ v_1\mapsto\zeta v_1,~ v_2\mapsto v_2,\\
\alpha\ :\ &u_1\mapsto u_2\mapsto (u_1u_2)^{-1},\\
&v_1\mapsto v_2\mapsto (v_1v_2)^{-1}.
\end{align*}}
From Theorem \ref{t2.2} it follows that if $K(u_i,v_i:1\leq i\leq
2)^{G}$ is rational over $K$, so is $K(x_0,y_0,u_i,v_i:1\leq i\leq
2)^{G}$ over $K$.

Define $w_1=v_1^3,w_2=v_2/v_1,V_1=w_2,V_2=(w_1w_2^2)^{-1}$. Then
$K(u_1,u_2,V_1,V_2)=K(u_1,u_2,v_1,v_2)^{\langle\alpha_1^3\rangle}$,
and the actions of $\alpha_1,\alpha_2$ and $\alpha$ on
$K(u_1,u_2,V_1,Vw_2)$ are
{\allowdisplaybreaks\begin{align*}\alpha_1\ :\
&u_1\mapsto u_1,~ u_2\mapsto\zeta^{-1} u_2,~ V_1\mapsto\zeta V_1,~ V_2\mapsto V_2,\\
\alpha_2\ :\
&u_1\mapsto\zeta^{-1} u_1,~ u_2\mapsto\zeta^{-1} u_2,~ V_1\mapsto\zeta^{-1} V_1,~ V_2\mapsto\zeta^{-1} V_2,\\
\alpha\ :\ &u_1\mapsto u_2\mapsto (u_1u_2)^{-1},\\
&V_1\mapsto V_2\mapsto (V_1V_2)^{-1}.
\end{align*}}
Define $U_1=u_1^3,U_2=u_2/u_1,W_1=V_1/u_1,W_2=V_2/u_2$. Then
$K(U_1,U_2,W_1,W_2)=K(u_1,u_2,V_1,V_2)^{\langle\alpha_2\rangle}$,
and the actions of $\alpha_1$ and $\alpha$ on $K(U_1,U_2,W_1,W_2)$
are {\allowdisplaybreaks\begin{align*}\alpha_1\ :\
&U_1\mapsto U_1,~ U_2\mapsto\zeta^{-1} U_2,~ W_1\mapsto\zeta W_1,~ W_2\mapsto\zeta W_2,\\
\alpha\ :\ &U_1\mapsto U_2^3U_1,~ U_2\mapsto (U_1U_2^2)^{-1},\\
&W_1\mapsto W_2\mapsto (W_1W_2)^{-1}.
\end{align*}}
Define $\tilde u_1=U_2,\tilde u_2=(U_1U_2^2)^{-1},\tilde
v_1=W_1\tilde u_1,\tilde v_2=W_2\tilde u_2$. Then $K(\tilde
u_1,\tilde u_2,\tilde v_1,\tilde v_2)=K(U_1,U_2,W_1,W_2)$, and the
actions of $\alpha_1$ and $\alpha$ on $K(\tilde u_1,\tilde
u_2,\tilde v_1,\tilde v_2)$ are
{\allowdisplaybreaks\begin{align*}\alpha_1\ :\
&\tilde u_1\mapsto\zeta^{-1} \tilde u_1,~ \tilde u_2\mapsto\zeta^{-1} \tilde u_2,~ \tilde v_1\mapsto\tilde v_1,~ \tilde v_2\mapsto\tilde v_2,\\
\alpha\ :\ &\tilde u_1\mapsto \tilde u_2\mapsto (\tilde u_1\tilde u_2)^{-1},\\
&\tilde v_1\mapsto \tilde v_2\mapsto (\tilde v_1\tilde v_2)^{-1}.
\end{align*}}
Define $\tilde U_1=\tilde u_1^3, \tilde U_2=\tilde u_2/\tilde u_1$.
Then $K(\tilde U_1,\tilde U_2,\tilde v_1,\tilde v_2)=K(\tilde
u_1,\tilde u_2,\tilde v_1,\tilde v_2)^{\langle\alpha_1\rangle}$, and
the action of $\alpha$ on $K(\tilde U_1,\tilde U_2,\tilde v_1,\tilde
v_2)$ is {\allowdisplaybreaks\begin{align*} \alpha\ :\ &\tilde
U_1\mapsto \tilde U_2^3\tilde U_1,~ \tilde U_2\mapsto (\tilde
U_1\tilde
U_2^2)^{-1},\\
&\tilde v_1\mapsto \tilde v_2\mapsto (\tilde v_1\tilde v_2)^{-1}.
\end{align*}}
Finally, define $\tilde w_1=\tilde U_2,\tilde w_2=(\tilde U_1\tilde
U_2^2)^{-1}$. Hence $\alpha(\tilde w_2)=(\tilde w_1\tilde
w_2)^{-1}$. According to Lemma \ref{l2.7} the action of $\alpha$ on
$K(\tilde w_1,\tilde w_2,\tilde v_1,\tilde v_2)$ can be linearized.


\medskip


\begin{thebibliography}{AAAA}
\bibitem[AHK]{AHK}
H. Ahmad, S. Hajja and M. Kang, Rationality of some projective
linear actions, {\it J. Algebra} {\bf 228} (2000), 643--658.
\bibitem[Be1]{Be1} H. A. Bender, A determination of the groups of order
$p^5$, {\it Ann. Math.}, {\bf 29} No. 1/4 (1927-1928), 61--72.
\bibitem[Be2]{Be2}
H. A. Bender, On groups of order $p^m,p$ being an odd prime number,
which contain an abelian subgroup of order $p^{m-1}$, {\it Ann.
Math.}, {\bf 29} No. 1/4 (1927-1928), 88--94.
\bibitem[Bo]{Bo}
F. A. Bogomolov, The Brauer group of quotient spaces by linear group
actions, {\it Math. USSR Izv.} {\bf 30} (1988), 455–-485.
\bibitem[CK]{CK}
H. Chu and M. Kang, Rationality of $p$-group actions, {\it J.
Algebra} {\bf 237} (2001), 673--690.
\bibitem[GMS]{GMS}
S. Garibaldi, A. Merkurjev and J-P. Serre, Cohomological invariants
in Galois cohomology, AMS Univ. Lecture Series vol. 28, Amer. Math.
Soc., Providence, 2003.
\bibitem[HK]{HK}
S. Hajja and M. Kang, Some actions of symmetric groups, {\it J.
Algebra} {\bf 177} (1995), 511--535.
\bibitem[HKK]{HKK}
A. Hoshi, M. Kang, B. E. Kunyavskii, Noether's problem and
unramified Brauer groups, arXiv:1202.5812.
\bibitem[HuK]{HuK}
S. J. Hu and M. Kang, Noether's problem for some $p$-groups, in
``Cohomological and geometric approaches to rationality problems",
edited by F. Bogomolov and Y. Tschinkel, Progress in Math. vol. 282,
Birkh\"auser, Boston, 2010.
\bibitem[Ja]{Ja}
R. James, The groups of order $p^6$ ($p$ an odd prime), {\it Math.
Comp.} {\bf 34} No. 150 (1980), 613--637.
\bibitem[Ka1]{Ka1}
M. Kang, Noether's problem for metacyclic $p$-groups, {\it Adv.
Math.} {\bf 203} (2005), 554--567.
\bibitem[Ka2]{Ka2}
M. Kang, Noether's problem for $p$-groups with a cyclic subgroup of
index $p^2$, {\it Adv. Math.} {\bf 226} (2011) 218--234.
\bibitem[Ka3]{Ka3}
M. Kang, Rationality problem for some meta-abelian groups, {\it J.
Algebra} {\bf 322} (2009), 1214-1219.
\bibitem[KP]{KP}
M. Kang and B. Plans, Reduction theorems for Noether's problem, {\it
Proc. Amer. Math. Soc.} {\bf 137} (2009), 1867--1874.
\bibitem[MM]{MM}
J.M. Masley, H.L. Montgomery, Cyclotomic fields with unique
factorization, {\it J. Reine Angew. Math.} {\bf 286/287} (1976)
248-–256.
\bibitem[Sa1]{Sa1}
D. J. Saltman, Generic Galois extensions and problems in field
theory, {\it Adv. Math.} {\bf 43} (1982), 250--283.
\bibitem[Sa2]{Sa2}
D. J. Saltman, Noether's problem over an algebraically closed field,
{\it Invent. Math.} {\bf 77}  (1984), 71--84.
\bibitem[Sw]{Sw}
R. Swan, Noether's problem in Galois theory, in "Emmy Noether in
Bryn Mawr", edited by B. Srinivasan and J. Sally, Springer-Verlag,
Berlin, 1983.
\end{thebibliography}
\end{document}